\documentclass[12pt]{article}
\usepackage{graphicx}
\usepackage{amsmath}
\usepackage[dvips]{color}
\usepackage{amsfonts}
\usepackage{amssymb}
\newtheorem{theorem}{Theorem}

\newtheorem{corollary}[theorem]{Corollary}

\newtheorem{lemma}[theorem]{Lemma}

\newtheorem{proposition}[theorem]{Proposition}

\begin{document}

\title{Conformal equivalent metrics in a plane domain are uniquely determined from
travel times}
\author{Victor .Palamodov}
\date{}
\maketitle

\section{Travel time and hodograph}

Let $E$ be an Euclidean space with the line element $\mathrm{d}s$ and
$\Omega\subset E$ be a closured bounded domain with $C^{1}$-boundary $\Gamma.$
Let $\mathbf{n}$ be a positive $C^{1}$-function in the closure $\Omega$
(refraction coefficient); consider the conformal metric $\mathbf{g}%
=\mathbf{n}^{2}\mathrm{d}s^{2}$ in $\Omega.$ Take some points $x,y\in\Omega$
and consider a geodesic curve $\gamma$ between these points. The length
$\tau\left(  \gamma\right)  $ of $\gamma$ is called travel time. It is a
multi-valued function $\tau\left(  x,y\right)  $ of the end points $x$ and $y$
of $\gamma$, if, at least one geodesic has a pair of conjugated points. This
is typical situation, if $\mathbf{n}$ is not constant, moreover, conjugated
points appear for any $\mathbf{n}\neq\mathrm{const}$, if $\mathbf{n}$ is
constant at infinity, see \cite{Cr}. The restriction $T$ of the travel time
function to $\Gamma\times\Gamma$ is called the \textit{hodograph. }The
branches of the hodograph $T\left(  x,y\right)  $ can be distinguished by
values of the tangent vector $\theta$ at one end, say $x.$ This vector can be
found from the differential of the corresponding branch of the hodograph, see Sec.7.

\textbf{Problem: }whether the refraction coefficient $\mathbf{n}$ in $\Omega$
\ is uniquely determined from knowledge of its hodograph $T.$

In other words, if two conformal metrics in $\Omega$ have equal hodographs,
does it imply that the metrics coincide?

This problem was studied by Muhometov \cite{Mu1}, Muhometov and Romanov
\cite{MR}, Beylkin \cite{Be2}, Bernstein and Gerver \cite{BG1}\cite{BG2}. The
stability of reconstruction of the metric from hodograph was proved in these
papers under assumption that the travel time $\tau$ is uniquely defined for
any pair of points $x,y\in\bar{\Omega},$ that is the hodograph is a usual
function. Then the equation $\tau_{1}\left(  x,y\right)  =\tau_{2}\left(
x,y\right)  $ for all pairs of points on the boundary implies $\mathbf{n}%
_{1}=\mathbf{n}_{2}$ in $\Omega.$ See also the papers of Croke \cite{Cr},
Sylvester and Uhlmann \cite{SyU}. More references are given in the surveys
\cite{Ger} and \cite{SU}.

We prove here an estimate similar to that of \cite{Mu1} $\left(  n=2\right)  $
without assumption of uniqueness of the travel time.

\section{Conformal metric}

Take the problem in a more general setting. Let $\Omega$ be a closed bounded
domain in $E=\mathbb{R}^{2}$ with the boundary $\Gamma$ of class $C^{2}.$ Let
$g=g_{ij}$\textrm{d}$x^{i}$\textrm{d}$x^{j}$ be a Riemannian metric in
$\Omega$ of class $C^{2};$ $x^{1},x^{2}$ are linear coordinates in $E.$ Let
$\mathbf{n\in} C^{2}\left(  \Omega\right)  $ be a positive function, called
refraction coefficient; consider the conformal Riemannian metric $\mathbf{g=n}g.$

\textbf{Assumptions }on a refraction coefficient: (i) $\mathbf{n}\in
C^{2}\left(  \Omega\right)  ,$\newline \indent(ii) for any point $x\in\Omega$
any geodesics $\gamma$ that starts at a point $x\in\Omega$ reaches the
boundary $\Gamma$ with a non tangent direction. This implies that there is no
waveguide in $\Omega$ (non-trapping geometry).

Let $T\left(  E\right)  $ be the tangent bundle on $E,$ $S\left(
\Omega\right)  $ be the bundle of unit circles in $T\left(  \Omega\right)  $
with respect to $g.$ Let $\gamma=\gamma\left(  x,\theta\right)  $ be the
geodesic curve that starts from a point $y=y\left(  x,\theta\right)  \in
\Gamma$ and arrives at the point $x$ with the unit tangent vector $\theta$.
The parameterization $\gamma=\gamma\left(  x,\theta\right)  ,\left(
x,\theta\right)  \in S\left(  \Omega\right)  $ can be applied to all geodesics
in $\Omega$ and the function
\begin{equation}
\tau\left(  x,\theta\right)  =\int_{\gamma\left(  x,\theta\right)  }%
\mathbf{n}\mathrm{d}s \label{tn}%
\end{equation}
is defined and is $C^{1}$-continuous on $S\left(  \Omega\right)  .$ It
vanishes for outgoing geodesics, that is for $\left\langle \nu|\theta
\right\rangle \leq0,$ where $\nu\left(  x\right)  $ means the inward conormal
to $\Gamma$ at $x.$ For incoming geodesic we have $\left\langle \nu
|\theta\right\rangle >0.$ So the data of travel times $\tau\left(
x,\theta\right)  \,$for $x\in\Gamma$ is the hodograph of the metric
$\mathbf{g}$.

\section{Differential of the travel time}

Let $\mathbf{n}$ be a refraction coefficient in $\Omega$ that fulfil the above
assumptions (i) and (ii). Fix an arbitrary point $\left(  x,\theta\right)  \in
S\left(  \Omega\right)  ;$ let $\gamma\left(  x,\theta\right)  $ be the
geodesic of the metric $\mathbf{g=n}^{2}g$ that arrive at $x$ with the
$g$-unit tangent vector $\theta$ and $y\in\Gamma$ be its initial point. Fix a
vector $\theta$ and consider $\tau\left(  x,\theta\right)  $ as a function of
$x$ only.

\begin{proposition}
\label{D}For any point $\left(  x,\theta\right)  \in S\left(  \Omega\right)
\,$ we have
\begin{equation}
\mathrm{d}_{x}\tau\left(  x,\theta\right)  =\mathbf{n}\left(  x\right)
\frac{\mathbf{e}\left(  x,\theta\right)  }{\left\langle \mathbf{e}\left(
x,\theta\right)  |\theta\right\rangle }. \label{dt}%
\end{equation}
where $\mathbf{e}\left(  x,\theta\right)  $ is a unit covector such that
$\left\langle \mathbf{e}\left(  x,\theta\right)  |\theta\right\rangle >0.$
\end{proposition}

\textsc{Proof.} Write $\mathrm{d}_{x}\tau\left(  x,\theta\right)
=\lambda\mathbf{n}\left(  x\right)  \mathbf{e}\left(  x,\theta\right)  $ for a
unit covector $\mathbf{e}$ and a scalar $\lambda>0.$ We have $\left\langle
\mathrm{d}_{x}\tau\left(  x,\theta\right)  |\theta\right\rangle =\mathbf{n}%
\left(  x\right)  >0$ according to the (\ref{tn}), which implies
$\lambda=\left\langle \mathbf{e}\left(  x,\theta\right)  |\theta\right\rangle
^{-1}.$ $\blacktriangleright$

We have $0<\left\langle \mathbf{e}\left(  x,\theta\right)  |\theta
\right\rangle \leq1$ and can write $\left\langle \mathbf{e}\left(
x,\theta\right)  |\theta\right\rangle =$ $\cos\omega$ for some angle
$\omega\in\left(  -\pi/2,\pi/2\right)  .$

Choose an orientation in $E;$ it induces an orientation in the tangent bundle
$T\left(  E\right)  .$ Take an orthonormal basis $\left(  a,b\right)  $ in the
cotangent bundle $T^{\ast}\left(  \Omega\right)  $ that is consistent with
this orientation. For a tangent vector $\theta$ we denote by $\hat{\theta}$
the dual covector; its coordinate expression is $\hat{\theta}_{i}=g_{ij}%
\theta^{j}.$ Write this covector in the form $\hat{\theta}=\cos\phi
\,a+\sin\phi\,b$ for some angle $\phi\in S^{1},$ which is a coordinate in the
bundle of all unit covectors in $\Omega.$ The coordinate system $x^{1}%
,x^{2},\phi$ is globally defined in this bundle. Let $\mathrm{d}$ be the
exterior differential in the complex of differential forms on $S^{\ast}\left(
\Omega\right)  ;$ write $\mathrm{d}=\mathrm{d}_{x}+\mathrm{d}_{\phi},$ where
$\mathrm{d}_{\phi}f=f^{\prime}\mathrm{d}\phi,$ $f^{\prime}\doteq\partial
f/\partial\phi.$ Differentiating (\ref{dt}) yields
\begin{equation}
\mathrm{dd}_{x}\tau=\mathrm{d}_{\phi}\mathrm{d}_{x}\tau=\mathbf{n}\left(
x\right)  \left[  \frac{\mathrm{d}_{\phi}\mathbf{e}\left(  x,\theta\right)
}{\cos\omega}+\frac{\sin\omega}{\cos^{2}\omega}\mathrm{d}_{\phi}\omega
\wedge\mathbf{e}\left(  x,\theta\right)  \right]  \label{ddt}%
\end{equation}
Set $\hat{\eta}=-\sin\phi\,a+\cos\phi\,b.$ The covectors $\hat{\theta}%
,\hat{\eta}$ form a positively oriented orthogonal frame and satisfy the
equations
\begin{align*}
\hat{\theta}\wedge\hat{\eta}  &  =\mathrm{d}V,\;\mathrm{d}\hat{\theta}%
=\hat{\eta}\wedge\mathrm{d}\phi,\\
\hat{\eta}\wedge\mathrm{d}\hat{\eta}  &  =\hat{\theta}\wedge\mathrm{d}%
\hat{\theta}=\mathrm{d}\phi\mathrm{d}V,
\end{align*}
where $\mathrm{d}V=\left(  \det\left\{  g_{ij}\right\}  \right)
^{1/2}\mathrm{d}x^{1}\mathrm{d}x^{2}$ is the Riemannian volume form in
$\left(  \Omega,g\right)  .$ Fix the angle $\omega$ by the equation
$\mathbf{e}=\cos\omega\,\hat{\theta}+\sin\omega\,\hat{\eta}.$

\section{Comparing two hodographs}

Now we estimate difference between two positive functions $\mathbf{n}%
_{1},\mathbf{n}_{2}$ (refraction coefficients) in terms of hodographs of two
conformal metrics $\mathbf{g}_{1}=\mathbf{n}_{1}^{2}g,$ $\mathbf{g}%
_{2}=\mathbf{n}_{2}^{2}g$. Denote by $\gamma_{1}\left(  x,\theta\right)
,\gamma_{2}\left(  x,\theta\right)  $ the corresponding geodesics, by
$\tau_{1}\left(  x,\theta\right)  $, $\tau_{2}\left(  x,\theta\right)  $ the
corresponding travel times. By Proposition \ref{D} we have $\mathrm{d}_{x}%
\tau_{j}\left(  x,\theta\right)  =\mathbf{n}_{j}\left(  x\right)  \cos
^{-1}\omega_{j}\mathbf{e}_{j}\left(  x,\theta\right)  $ for some angles
$\omega_{1},\omega_{2}$ such that $\mathbf{e}_{j}=\cos\omega_{j}\hat{\theta
}+\sin\omega_{j}\hat{\eta}.$ Set $\rho\left(  x,\phi\right)  =\tau_{2}\left(
x,\theta\right)  -\tau_{1}\left(  x,\theta\right)  $ and calculate the product%
\[
\mathbf{R}\doteq\mathrm{d}\rho\wedge\mathrm{d}_{\phi}\mathrm{d}\rho=\left(
\mathrm{d}\tau_{2}-\mathrm{d}\tau_{1}\right)  \wedge\left(  \mathrm{d}_{\phi
}\mathrm{d}\tau_{2}-\mathrm{d}_{\phi}\mathrm{d}\tau_{1}\right)
\]

\begin{lemma}
\label{l}We have
\begin{align}
\mathbf{R}  &  =\left[  \left(  \frac{\mathbf{n}_{2}}{\cos\omega_{2}}\right)
^{2}+\left(  \frac{\mathbf{n}_{1}}{\cos\omega_{1}}\right)  ^{2}-\frac
{2\cos\left(  \omega_{1}-\omega_{2}\right)  \mathbf{n}_{1}\mathbf{n}_{2}}%
{\cos\omega_{1}\cos\omega_{2}}\right]  \hat{\theta}\wedge\mathrm{d}\hat
{\theta}\label{rnn}\\
&  +\mathrm{d}_{\phi}\left(  \mathbf{n}_{2}\tan\omega_{2}-\mathbf{n}_{1}%
\tan\omega_{1}\right)  \left(  \mathbf{n}_{2}-\mathbf{n}_{1}\right)
\wedge\hat{\theta}\wedge\hat{\eta}\nonumber
\end{align}
\end{lemma}

\textsc{Proof. }Substituting (\ref{ddt}) yields%
\begin{align*}
\mathbf{R}  &  =\mathbf{n}_{2}^{2}\frac{\mathbf{e}_{2}\wedge\mathrm{d}%
\mathbf{e}_{2}}{\cos^{2}\omega_{2}}+\mathbf{n}_{1}^{2}\frac{\mathbf{e}%
_{1}\wedge\mathrm{d}\mathbf{e}_{1}}{\cos^{2}\omega_{1}}-\mathbf{n}%
_{1}\mathbf{n}_{2}\frac{\mathbf{e}_{2}\wedge\mathrm{d}\mathbf{e}%
_{1}+\mathbf{e}_{1}\wedge\mathrm{d}\mathbf{e}_{2}}{\cos\omega_{1}\cos
\omega_{2}}\\
&  -\mathbf{n}_{1}\mathbf{n}_{2}\left[  \frac{\sin\omega_{1}\mathrm{d}%
\omega_{1}}{\cos^{2}\omega_{1}\cos\omega_{2}}+\frac{\sin\omega_{2}%
\mathrm{d}\omega_{2}}{\cos\omega_{1}\cos^{2}\omega_{2}}\right]  \wedge
\mathbf{e}_{1}\wedge\mathbf{e}_{2}.
\end{align*}
We have for $\mathbf{e=e}_{j},j=1,2$
\begin{align*}
\mathbf{e}\wedge\mathrm{d}\mathbf{e}  &  =\left(  \cos\omega\hat{\theta}%
+\sin\omega\hat{\eta}\right)  \wedge\mathrm{d}\left(  \cos\omega\hat{\theta
}+\sin\omega\hat{\eta}\right) \\
&  =\left(  \cos\omega\hat{\theta}+\sin\omega\hat{\eta}\right)  \wedge\left(
\cos\omega\mathrm{d}\hat{\theta}+\sin\omega\mathrm{d}\hat{\eta}-\sin
\omega\mathrm{d}\omega\hat{\theta}+\cos\omega\mathrm{d}\omega\hat{\eta}\right)
\\
&  =\cos^{2}\omega\hat{\theta}\wedge\mathrm{d}\hat{\theta}+\sin^{2}\omega
\hat{\eta}\wedge\mathrm{d}\hat{\eta}+\hat{\theta}\wedge\hat{\eta}%
\mathrm{d}\omega\\
&  =\hat{\theta}\wedge\mathrm{d}\hat{\theta}+\hat{\theta}\wedge\hat{\eta
}\mathrm{d}\omega=\left(  1+\omega^{\prime}\right)  \hat{\theta}%
\wedge\mathrm{d}\hat{\theta},
\end{align*}
since $\hat{\theta}\wedge\mathrm{d}\hat{\eta}=\mathrm{d}\hat{\theta}\wedge
\hat{\eta}=0,$ and
\begin{align*}
\mathbf{e}_{1}\wedge\mathrm{d}\mathbf{e}_{2}  &  =\left(  \cos\omega_{1}%
\hat{\theta}+\sin\omega_{1}\hat{\eta}\right)  \wedge\mathrm{d}\left(
\cos\omega_{2}\hat{\theta}+\sin\omega_{2}\hat{\eta}\right) \\
&  =\cos\left(  \omega_{1}-\omega_{2}\right)  \hat{\theta}\wedge\mathrm{d}%
\hat{\theta}+\cos\left(  \omega_{1}-\omega_{2}\right)  \hat{\theta}\wedge
\hat{\eta}\mathrm{d}\omega_{2}\\
&  =\cos\left(  \omega_{1}-\omega_{2}\right)  \left(  1+\omega_{2}^{\prime
}\right)  \hat{\theta}\wedge\mathrm{d}\hat{\theta},\\
\mathbf{e}_{1}\wedge\mathbf{e}_{2}  &  =-\sin\left(  \omega_{1}-\omega
_{2}\right)  \mathrm{d}V.
\end{align*}
Therefore
\begin{align*}
\frac{\mathbf{R}}{\hat{\theta}\wedge\mathrm{d}\hat{\theta}}  &  =\mathbf{n}%
_{2}^{2}\frac{1+\omega_{2}^{\prime}}{\cos^{2}\omega_{2}}+\mathbf{n}_{1}%
^{2}\frac{1+\omega_{1}^{\prime}}{\cos^{2}\omega_{1}}-\mathbf{n}_{1}%
\mathbf{n}_{2}\frac{\cos\left(  \omega_{1}-\omega_{2}\right)  \left(
2+\omega_{2}^{\prime}+\omega_{1}^{\prime}\right)  }{\cos\omega_{1}\cos
\omega_{2}}\\
&  -\mathbf{n}_{1}\mathbf{n}_{2}\sin\left(  \omega_{2}-\omega_{1}\right)
\frac{\cos\omega_{2}\sin\omega_{1}\omega_{1}^{\prime}-\cos\omega_{1}\sin
\omega_{2}\omega_{2}^{\prime}}{\cos^{2}\omega_{2}\cos^{2}\omega_{1}}.
\end{align*}
Two terms containing $\phi$-derivatives gives%

\begin{align*}
&  \frac{\cos\left(  \omega_{1}-\omega_{2}\right)  \left(  \omega_{2}^{\prime
}+\omega_{1}^{\prime}\right)  }{\cos\omega_{1}\cos\omega_{2}}+\sin\left(
\omega_{2}-\omega_{1}\right)  \frac{\cos\omega_{2}\sin\omega_{1}\omega
_{1}^{\prime}-\cos\omega_{1}\sin\omega_{2}\omega_{2}^{\prime}}{\cos^{2}%
\omega_{2}\cos^{2}\omega_{1}}\\
&  =\frac{\omega_{1}^{\prime}}{\cos^{2}\omega_{1}}+\frac{\omega_{2}^{\prime}%
}{\cos^{2}\omega_{2}},
\end{align*}
times $\mathbf{n}_{1}\mathbf{n}_{2},$ which finally yields
\begin{align*}
\frac{\mathbf{R}}{\hat{\theta}\wedge\mathrm{d}\hat{\theta}}  & =\left(
\frac{\mathbf{n}_{2}}{\cos\omega_{2}}\right)  ^{2}+\left(  \frac
{\mathbf{n}_{1}}{\cos\omega_{1}}\right)  ^{2}-\frac{2\cos\left(  \omega
_{1}-\omega_{2}\right)  \mathbf{n}_{1}\mathbf{n}_{2}}{\cos\omega_{1}\cos
\omega_{2}}\\
&  +\frac{\mathbf{n}_{2}^{2}\omega_{2}^{\prime}}{\cos^{2}\omega_{2}}%
+\frac{\mathbf{n}_{1}^{2}\omega_{1}^{\prime}}{\cos^{2}\omega_{1}}%
-\mathbf{n}_{1}\mathbf{n}_{2}\left[  \frac{\omega_{1}^{\prime}}{\cos^{2}%
\omega_{1}}+\frac{\omega_{2}^{\prime}}{\cos^{2}\omega_{2}}\right]  ,
\end{align*}
This completes the proof. $\blacktriangleright$

\section{The basic inequality}

\begin{theorem}
\label{main}If $\mathbf{n}_{1}$, $\mathbf{n}_{2}$ fulfil the conditions (i)
and (ii) and $\mathbf{n}_{1}=\mathbf{n}_{2}$ on $\Gamma,$ the inequality
holds
\begin{equation}
\int\left(  \frac{\mathbf{n}_{2}}{\cos\omega_{2}}-\frac{\mathbf{n}_{1}}%
{\cos\omega_{1}}\right)  ^{2}\mathrm{d}\phi\mathrm{d}V\leq-\int_{S^{1}}%
\int_{\Gamma}\mathrm{d}_{x}\rho\wedge\mathrm{d}_{\phi}\rho. \label{ineq}%
\end{equation}
\end{theorem}

\textsc{Proof. }The last term of (\ref{rnn}) is equal to the exact form
\[
\mathrm{d}\left[  \left(  \mathbf{n}_{2}\tan\omega_{2}-\mathbf{n}_{1}%
\tan\omega_{1}\right)  \left(  \mathbf{n}_{2}-\mathbf{n}_{1}\right)
\hat{\theta}\wedge\hat{\eta}\right]  ,
\]
since $\mathrm{d}\left(  \hat{\theta}\wedge\hat{\eta}\right)  =\mathrm{d}%
\left(  a\wedge b\right)  =0.$ It is continuous up the boundary $\Gamma$ as
well as its primitive function $\left(  \mathbf{n}_{2}\tan\omega
_{2}-\mathbf{n}_{1}\tan\omega_{1}\right)  \left(  \mathbf{n}_{2}%
-\mathbf{n}_{1}\right)  $ in spite of the functions $\tan\omega_{1},\tan
\omega_{2}$ tend to infinity at $\Gamma$. This follows from vanishing of
$\mathbf{n}_{1}-\mathbf{n}_{2}$ on $\Gamma.$ Therefore integration over
$S\left(  \Omega\right)  $ gives
\begin{align*}
& \int_{S\left(  \Omega\right)  }\mathrm{d}\rho\wedge\mathrm{d}_{\phi
}\mathrm{d}\rho\\
= & \int_{S\left(  \Gamma\right)  }\left[  \left(  \frac{\mathbf{n}_{2}}%
{\cos\omega_{2}}\right)  ^{2}+\left(  \frac{\mathbf{n}_{1}}{\cos\omega_{1}%
}\right)  ^{2}-\frac{2\cos\left(  \omega_{1}-\omega_{2}\right)  \mathbf{n}%
_{1}\mathbf{n}_{2}}{\cos\omega_{1}\cos\omega_{2}}\right]  \mathrm{d}%
\phi\mathrm{d}V\\
&  \geq\int_{S\left(  \Gamma\right)  }\left[  \left(  \frac{\mathbf{n}_{2}%
}{\cos\omega_{2}}\right)  ^{2}+\left(  \frac{\mathbf{n}_{1}}{\cos\omega_{1}%
}\right)  ^{2}-\frac{2\mathbf{n}_{1}\mathbf{n}_{2}}{\cos\omega_{1}\cos
\omega_{2}}\right]  \mathrm{d}\phi\mathrm{d}V\\
&  =\int_{S\left(  \Gamma\right)  }\left(  \frac{\mathbf{n}_{2}}{\cos
\omega_{2}}-\frac{\mathbf{n}_{1}}{\cos\omega_{1}}\right)  ^{2}\mathrm{d}%
\phi\mathrm{d}V
\end{align*}
By Stokes'
\[
\int_{S\left(  \Omega\right)  }\mathrm{d}\rho\wedge\mathrm{d}_{\phi}%
\mathrm{d}\rho=-\int_{S\left(  \Omega\right)  }\mathrm{d}\left(
\mathrm{d}_{x}\rho\wedge\mathrm{d}_{\phi}\rho\right)  =-\int_{S\left(
\Gamma\right)  }\mathrm{d}_{x}\rho\wedge\mathrm{d}_{\phi}\rho,
\]
since $\mathrm{d}_{\phi}\mathrm{d}=-\mathrm{dd}_{x},$ and (\ref{ineq})
follows. $\blacktriangleright$

\begin{corollary}
If the hodographs of $\mathbf{n}_{1}$ and $\mathbf{n}_{2}$ coincide for all
incoming directions $\left(  x,\theta\right)  ,x\in\Gamma,$ then
$\mathbf{n}_{2}=\mathbf{n}_{1}.$
\end{corollary}

\textsc{Proof of Corollary.} It is easy to see that the condition implies
equation $\mathbf{n}_{1}=\mathbf{n}_{2}$ on $\Gamma.$ Therefore we can apply
Theorem \ref{main} and conclude that $\rho\left(  x,\phi\right)  \equiv0.$
This yields $\mathbf{n}_{2}\cos\omega_{1}\equiv\mathbf{n}_{1}\cos\omega_{2}$
for all $x\in\Omega$ and unit vectors $\theta.$ Fix $x$ and choose $\theta$ in
such a way that the geodesic $\gamma_{1}\left(  x,\theta\right)  $ has maximal
length. The geodesic $\gamma_{1}$ is orthogonal to $\Gamma$ at the initial
point $y_{1}$ and by the Gauss Lemma (see \cite{BK}) the end points $z$ of
geodesics $\gamma_{1}\left(  z,\theta\right)  $ such that $\tau_{1}\left(
z,\theta\right)  =\tau_{1}\left(  x,\theta\right)  $ run over a curve $I$ that
is orthogonal to $\gamma_{1}\left(  x,\theta\right)  .$ Therefore
$\mathbf{e}_{1}=\hat{\theta}$ and $\omega_{1}=0.$ The equation$\,\cos
\omega_{1}=1$ implies $\mathbf{n}_{2}\left(  x\right)  \leq\mathbf{n}%
_{1}\left(  x\right)  .$ The opposite inequality is also true.
$\blacktriangleright$

\section{Application to the geodesic transform}

\begin{corollary}
For an arbitrary metric $\mathbf{n}$ that fulfils (i) and (ii) and any real
$C^{1}$-smooth function $f$ in $\Omega$ that vanishes at $\Gamma$ we have
\begin{equation}
\int\frac{\mathrm{d}\phi}{\cos^{2}\omega}f^{2}\left(  x\right)  \mathrm{d}%
V\leq-\int_{S^{1}}\int_{\Gamma}\mathrm{d}_{x}g\wedge\mathrm{d}_{\phi}g
\label{fg}%
\end{equation}
where
\[
g\left(  x,\theta\right)  =\int_{\gamma\left(  x,\theta\right)  }%
f\mathrm{d}s.
\]
\end{corollary}

\textsc{Proof. }We may assume that $f\in C^{2}\left(  \Omega\right)  $ and
apply theorem \ref{main} to $\mathbf{n}_{1}=\mathbf{n,}$ $\mathbf{n}%
_{2}=\mathbf{n}+\varepsilon f,$ where $\varepsilon$ is a small parameter. The
geodesic curves for $\mathbf{n}_{1}$and $\mathbf{n}_{2}$ are the same up to
$O\left(  \varepsilon^{2}\right)  .$ Collecting the terms of order
$\varepsilon$ in (\ref{ineq}), we obtain (\ref{fg}). $\blacktriangleright$

\textbf{Example.} Consider the Euclidean metric $\mathbf{g}=\mathrm{d}s^{2}$
in the unit disc $\Omega$. By a direct calculation
\[
\int\frac{\mathrm{d}\phi}{\cos^{2}\omega}=2\pi\left(  \frac{1+\left|
x\right|  }{1-\left|  x\right|  }\right)  ^{1/2},\left|  x\right|  <1.
\]
Write the line integrals $g$ in the the standard parameterization $G\left(
p,\varphi\right)  =g\left(  x,\theta\right)  $ of line integrals, where
\begin{align*}
x  &  =\left(  \cos s,\sin s\right)  ,\,-\pi/2\leq s\leq\pi/2,\theta=\left(
\cos\phi,\sin\phi\right)  ,s-\pi/2<\phi<s+\pi/2,\\
p  &  =\cos\left(  \psi-\phi\right)  ,\varphi=\phi+\pi/2
\end{align*}
and
\[
g_{s}^{\prime}=-\sqrt{1-p^{2}}G_{p}^{\prime},\;g_{\phi}^{\prime}=\sqrt
{1-p^{2}}G_{p}^{\prime}+G_{\varphi}^{\prime}%
\]
The inequality (\ref{fg}) takes the form
\begin{align*}
\int\left(  \frac{1+\left|  x\right|  }{1-\left|  x\right|  }\right)
^{1/2}f^{2}\left(  x\right)  \mathrm{d}x  &  \leq-\int_{0}^{2\pi}%
\mathrm{d}s\int_{s-\pi/2}^{s+\pi/2}g_{s}^{\prime}\left(  x,\theta\right)
g_{\phi}^{\prime}\left(  x,\theta\right)  \mathrm{d}\phi\\
&  =\int_{0}^{2\pi}\int_{-1}^{1}\left(  \sqrt{1-p^{2}}\left(  G_{p}^{\prime
}\right)  ^{2}+G_{p}^{\prime}G_{\varphi}^{\prime}\right)  \mathrm{d}%
p\mathrm{d}\varphi.
\end{align*}

\section{Calculation of angular derivative of the travel time}

Take an arbitrary refraction coefficient $\mathbf{n}$ in $\Omega$ that
satisfies (i) and (ii) and determine the angular parameter $\theta$ of a ray
$\gamma$ in terms of the hodograph function. Let $s$ be the arc length along
$\Gamma$ that increases clockwise.

\begin{proposition}
Let $\left(  x,\theta\right)  ,x\in\Gamma$ be an arbitrary outward direction
and $\gamma\left(  s\right)  $ be the family of geodesic that joins some point
$y\in\Gamma$ and $x\left(  s\right)  \in\Gamma$ with outgoing tangent vector
$\theta\left(  s\right)  $ defined for $s\in\left[  0,\varepsilon\right]  $
(or for $s\in\left[  -\varepsilon,0\right]  $) such that $x\left(  0\right)
=x,$ $\theta\left(  0\right)  =\theta.$ Then we have
\[
\sin\psi=\frac{1}{\mathbf{n}\left(  x\right)  }\frac{\partial\tau\left(
\gamma\left(  s\right)  \right)  }{\partial s}|_{s=0},
\]
where $\psi$ is the angle between the outward conormal $-\nu\left(  x\right)
$ and the outward covector $\hat{\theta}.$
\end{proposition}

\textsc{Proof. }By the Gauss Lemma the tangent vector $\theta$ is orthogonal
to the circle of points with the $\mathbf{g}$-distance $\tau\left(
x,y\right)  $ from $y.$ $\blacktriangleright$

This formula also holds for any point $x$ conjugate to $y$ in $\gamma\left(
0\right)  .$ In particular, $\theta=-\nu,$ if $\partial\tau/\partial s=0,$
that is the geodesic $\gamma\left(  0\right)  $ arrives to $x$ with the normal
direction $-\nu\left(  x\right)  $. The angle $\phi$ of $\theta$ of the ray
$\gamma\left(  x,\theta\right)  $ is equal to $\phi=\psi+\beta+\pi,$ where
$\beta=\mathrm{\arg\,}\nu.$

\end{document}